\documentclass[11pt]{amsart}
\usepackage{latexsym}
\usepackage{amssymb}
\usepackage{amsxtra}
\usepackage{amsmath}
\usepackage{amsfonts}

\newcommand{\cleqn}{\setcounter{equation}{0}}
\newcommand{\clth}{\setcounter{theorem}{0}} 
\newcommand {\sectionnew}[1]{\section{#1}\cleqn\clth}
\newcommand{\beq}{\begin{equation}}
\newcommand{\eeq}{\end{equation}}
\newcommand{\beqa}{\begin{eqnarray}}
\newcommand{\eeqa}{\end{eqnarray}}
\newcommand{\beaa}{\begin{eqnarray*}}
\newcommand{\ben}{\begin{eqnarray*}}
\newcommand{\eaa}{\end{eqnarray*}}
\newcommand{\een}{\end{eqnarray*}}

\newcommand \nc {\newcommand}

\newtheorem{theorem}{Theorem}[section]
\newtheorem{lemma}[theorem]{Lemma}
\newtheorem{proposition}[theorem]{Proposition}
\newtheorem{corollary}[theorem]{Corollary}

\nc \thref{Theorem \ref}
\nc \leref{Lemma \ref}
\nc \prref{Proposition \ref}
\nc \coref{Corollary \ref}
\nc \deref{Definition \ref}
\nc \exref{Example \ref}
\nc \reref{Remark \ref}

\newcommand{\CC}{\mathcal{C}}

\renewcommand{\P}{\mathcal{P}}
\newcommand{\Q}{\mathcal{Q}}

\newcommand{\W}{\mathcal{W}}
\newcommand{\A}{\mathcal{A}}

\newcommand{\C}{\mathbb{C}}

\renewcommand{\L}{\mathcal{L}}

\newcommand{\Z}{\mathbb{Z}}

\newcommand{\q}{\mathbf{q}}

\newcommand{\x}{\mathbf{x}}
\newcommand{\y}{\mathbf{y}}

\def\res{\mathop{\rm Res}\nolimits}

\def\d{\partial}

\def\tensor{\otimes}

\def\({\left(}
\def\){\right)}
\def\[{\left[}
\def\]{\right]}
\def\<{\left\langle}
\def\>{\right\rangle}

\def\gl{\lambda}
\def\ge{\epsilon}
\def\ga{\alpha}
\def\gd{\delta}
\def\gb{\beta}
\begin{document}
\title{Hirota quadratic equations for the Extended Toda Hierarchy }
\author{
Todor E. Milanov}
\thanks{E-mail: milanov@math.berkeley.edu}
\date{}

\begin{abstract}
The Extended Toda Hierarchy (shortly ETH) was introduce by E. Getzler
\cite{Ge} and independently by Y. Zhang \cite{Z} in order to describe 
an integrable 
hierarchy which governs the Gromov--Witten invariants of $\C P^1$. 
The {\em Lax type} presentation of the ETH was given in
\cite{CDZ}. In this paper we give a description of 
the ETH in terms of {\em tau-functions} and Hirota Quadratic Equations
(known also as Hirota Bilinear Equations). 
A new feature here is that the Hirota equations are given in terms of 
vertex operators taking values in the algebra of differential operators 
on the affine line. 
\end{abstract}    

\maketitle

\setcounter{section}{0}
\sectionnew{Introduction}

We begin with a description of the KdV hierarchy of integrable systems
which will be used as a model for introducing tau-functions and Hirota
Quadratic Equations (shortly HQE) of the ETH. Both the KdV and the ETH
hierarchies will be presented in the form convinient for the applications
to Gromov--Witten theory. In particular, we work over the 
field $\C((\ge)),$ where $\ge$ corresponds to the genus expansion 
parameter.   

The KdV hierarchy can be described in the {\em Lax form}
\cite{GD} as a sequence of commuting flows on the space of {\em Lax operators}
\begin{equation} \notag 
L = \frac{\ge^2\d_x^2}{2}\ +\ u(x,\ge), \end{equation}
where $u=u_0(x)+u_1(x)\ge+u_2(x)\ge^2+...$ is a power $\ge$-series
with infinitely differentiable coefficients. By definition, the flows have
the form
\begin{equation} \notag \frac{\d}{\d q_n} L = \ge^{-1}
\left[ \left(\frac{(2L)^{n+1/2}}{(2n+1)!!}\right)_{+}, L \right], 
\ \ \ n=0,1,2,..., \end{equation}
where $B_{+}$ means the differential part $\sum_{m\geq k\geq 0} w_k(x)\d_x^k$
of a pseudo-differential operator 
$B=\sum_{m\geq k >-\infty} w_k(x)\d_x^k$, and $(2n+1)!!=1\cdot 3\cdot 
5 \cdots (2n+1)$.

Given a Lax operator $L$ there exists an integral operator  
$P=1+w_1\d_x^{-1}+w_2\d_x^{-2}+\ldots$, called \emph{dressing operator} 
such that $L=P(\ge^2\d_x^2/2)P^{-1}$. A Lax operator $L$ is a solution to the
KdV hierarchy if and only if there exists a dressing operator $P$ satisfying
the equation
\ben
\ge \frac{\d}{\d q_n} P=-\left(\frac{(2L)^{n+1/2}}{(2n+1)!!}\right)_{-} P.
\een
Such a dressing operator is called \emph{wave operator of the KdV hierarchy}. 
The wave operator is uniquely determined by $L$ up to multiplication
from the right by an operator of the form $1+a_1\d_x^{-1}+a_2\d_x^{-2}+\ldots$
where $a_i$ are constants independent of $x$ and the time variables 
$q_n,\ n\geq 0$.

An equivalent description of the KdV hierarchy can be given in terms
of tau-functions and Hirota quadratic equations (see 
{\em e.g.} \cite{Kac}). A non-vanishing function
$\tau (q_0,q_1,... ; \ge)$, is called 
\emph{a tau-function} of the KdV hierarchy if there exists a 
wave operator  
$P = 1 + w_1 \d^{-1}_x+w_2\d_x^{-2}+... $ , such that
\begin{equation} \notag 
1+\frac{w_1}{\sqrt{2\gl}}+\frac{w_2}{\sqrt{2\gl}^2}+... := 
\frac{\exp\left(- \sum_{n\geq 0} \frac{(2n-1)!!}{(2\gl)^{1/2+n}} 
\ge \frac{\d}{\d q_n} \right) \tau (q_0+x,q_1,...;\ge)}
{\tau (q_0+x,q_1,...;\ge)} 
\end{equation}      
The tau-function is uniquely determined by $P$ up to multiplication by a 
constant (see {\em e.g.} \cite{Kac} or \cite{PvM}).

According to \cite{Kac} $\tau$ is a tau-function of the
KdV hierarchy if and only if it satisfies the following Hirota equation. 
  Put $\d_n =\d/\d q_n$ and consider the Heisenberg Lie algebra spanned 
by the operators of differentiation $\ge \d_n$ and multiplication 
$q_m/\ge$. The following {\em vertex operators} $\Gamma^{\pm}$ 
represent the action on the Fock space of certain elements of 
the Heisenberg group:
\beq \notag \Gamma^{\pm} := \exp \left( \pm \sum_{n\geq 0} 
\frac{(2\gl)^{n+1/2}}{(2n+1)!!} \frac{q_n}{\ge} \right)
\exp \left( \mp \sum_{n\geq 0} \frac{(2n-1)!!}{(2\gl)^{1/2+n}} \ge \d_n 
\right) .\eeq
We remark that 
\beq \notag 
\frac{(2n-1)!!}{(2\gl)^{1/2+n}} = \left(-\frac{d}{d\gl}\right)^n 
\frac{1}{\sqrt{2\gl}}, \ \ \frac{(2\gl)^{1/2+n}}{(2n+1)!!} = 
\left(\frac{d}{d\gl}\right)^{-1-n} \frac{1}{\sqrt{2\gl}} .\eeq

According to \cite{G} the Hirota equation of the KdV hierarchy can be 
stated this way:
\beq \label{KdV}
\frac{d\gl}{\sqrt{\gl}} \left( \Gamma^{+}\otimes \Gamma^{-} - 
\Gamma^{-}\otimes\Gamma^{+} \right)\ (\tau \otimes \tau) \ \ 
\text{\em is regular in $\gl$.} \eeq
Here $\tau \otimes \tau$  means the function $\tau (\q'; \ge) \tau (\q'';\ge)$
of the two copies of the variable $\q = (q_0,q_1,...)$, and the vertex
operators in $\Gamma^{\pm}\otimes\Gamma^{\mp}$ preceding (respectively --- 
following) $\otimes$ act on $\q'$ (respectively --- on $\q''$). 
The expression in (\ref{KdV}) is in fact single-valued in $\gl $ near 
$\gl=\infty$. Passing to the variables $\x = (\q'+\q'')/2$ and 
$\y = (\q'-\q'')/2$ and using Taylor's formula one can expand (\ref{KdV}) 
into a power series in $\y$ with coefficients which are Laurent series 
in $\gl^{-1}$ (whose coefficients are polynomials in $\tau(\x)$ and its
partial derivatives). The regularity condition in (\ref{KdV}) means,
by definition,  that all the Laurent series in $\gl^{-1}$ are polynomials
in $\gl$.

\medskip  

We describe the Lax form of the ETH following \cite{CDZ} 
\footnote{In the papers \cite{Ge,Z} the ETH is defined implicitly 
via bihamiltonian recursion relations}.  
Introduce a \emph{Lax operator}
\beq
\L=\Lambda +u+Qe^{v}\Lambda^{-1},
\label{laxop}
\eeq
where 
$\Lambda=e^{\epsilon\d_x}$ is the shift operator, i.e.
$
\Lambda a(x)=\sum_{k\geq 0} a^{(k)}(x)\epsilon^k/k!,
$
$u$ and $v$ are formal series
\ben
&& u=  u_0(x)+ u_1(x)\epsilon +u_2(x)\epsilon^2+\ldots \\
&& v= \ \ \ \ \ \ \ \ \ \ v_1(x)\epsilon +v_2(x)\epsilon^2+\ldots 
\een
with coefficients which are infinitely differentiable 
functions of $x$,
and $Q$ is a nonzero constant.

Let $\L$ be a Lax operator. A pair of \emph{dressing operators}
consists of two operators
\beqa
&&
\P_L=1+w_1\Lambda^{-1}+w_2\Lambda^{-2}+\ldots ,  
\label{dressP}\\
&&
\P_R=\tilde{w_0}+\Lambda^{-1}\tilde{w_1}+\Lambda^{-2}\tilde{w_2}+
\ldots ,  
\label{dressQ}
\eeqa
such that $\L=\P_L\Lambda \P_L^{-1} = \(\P_R^{-1} \Lambda\P_R\)^\# $, where
$\#$ is an {\em antiinvolution} acting on the space of Laurent series in
$\Lambda$ by $x^\#=x$ and $\Lambda^\# = Q\Lambda^{-1}$. 
The pair is unique up to multiplying $\P_L$ from the right and $\P_R$
from the left 
by operators of the form respectively 
$1+ a_1\Lambda^{-1}+a_2\Lambda^{-2}+...$ and 
$\tilde{a}_0 + \tilde{a}_1\Lambda +\tilde{a}_2\Lambda^2+\ldots$
with coefficients independent of $x$. 
Using $\P_L$ and $\P_R$ one expresses the logarithm of the Lax operator 
\ben
&
\log \L& = 
\frac{1}{2}(\P_L\epsilon \d_x\P_L^{-1} + \(\P_R^{-1}\epsilon\d_x\P_R\)^\#) =\\
&      & =
\frac{1}{2} \(\log Q +\epsilon\(\P_R^{-1}\frac{\d \P_R}{\d x}\)^\#-
                     \epsilon\frac{\d \P_L}{\d x}\P_L^{-1} \),
\een
as a Laurent series in $\Lambda$ possibly infinite in both directions.

By definition, the ETH is the following sequence of flows 
with time variables $q_{n,0}$ and $q_{n,1},\, n=0,1,2\ldots$:
\beqa
{\d_{n,1}}\L & = & 
\epsilon^{-1}
\left[ \left(  \frac{  \L^{n+1}  }{(n+1)!}  \right)_+,\L\right]
\label{laxeq1}                                               ,\\
{\d _{n,0}}\L    & = & 
2\epsilon^{-1}
\left[ \left(  
\frac{\L^n}{n!}(\log \L-{\CC}_n) \right)_+ ,\L \right].
\label{laxeq2}
\eeqa
where $\d_{n,i}:=\d/\d q_{n,i}$, ${\CC}_n$ are the harmonic numbers defined 
by: $\CC_0=(1/2)\log Q,$ $\CC_n= \CC_{n-1}+1/n$, 
and for $B=\sum B_k\Lambda^{k}$  we put 
$B_+=\sum_{k\geq 0}B_k\Lambda^k$.
According to \cite{CDZ}, the flows of ETH commute pairwise
and preserve the class of Lax operators \eqref{laxop} with a fixed $Q$.

We will prove that $\L$ is a solution to the ETH if and only if there is 
a pair of dressing operators $\P_L$ and $\P_R$, which satisfies the 
differential equations: 
{\allowdisplaybreaks 
\begin{align}
\label{bn1}
\d_{n,1}\P_L & =-  \left(  \frac{  \L^{n+1}  }{\ge(n+1)!}  \right)_- \P_L, \\
\label{bn0}
\d_{n,0}\P_L & =-\left(  
\frac{2\L^n}{\ge n!}(\log {\L}-{\CC}_n) \right)_-\P_L, \\
\label{an1} 
\d_{n,1}\P_R^\# & 
= \left(  \frac{  \L^{n+1}  }{\ge(n+1)!}  \right)_+\P_R^\# , \\
\label{an0}
\d_{n,0}\P_R^\# & = \left(  
\frac{2\L^n}{\ge n!}(\log {\L}-{\CC}_n) \right)_+ \P_R^\#.
\end{align}}
We call such a pair \emph{wave operators of the ETH}. It is unique up to
multiplying $\P_L$ from the right and $\P_R$ from the left  by operators
of the form respectively $1+a_1\Lambda^{-1}+a_2\Lambda^{-2}+\ldots$
and $\tilde a_0 +\tilde a_1\Lambda^{-1}+\tilde a_2\Lambda^{-2}+\ldots$ , 
where $a_i$ and $\tilde a_j$ are independent of $x$ and $\q$.

Given a non-vanishing function $\tau(\q;\ge)$ associate to it a pair of 
operators
\ben 
&&
\P_L=1+w_1\Lambda^{-1}+w_2\Lambda^{-2}+\dots ,\\
&&
\P_R=\tilde w_0 + \Lambda^{-1}\tilde w_1 + \Lambda^{-2} \tilde w_2 +\ldots 
\een
defined by
\ben
&&
P_L:=
1+\frac{w_1}{\gl}+\frac{w_2}{\gl^2}+\ldots :=  
\frac{\exp\(-\sum_{n\geq 0} n!\lambda^{-n-1} \ge\d_{n,1}\)
\tau (q_{0,0}+x-\frac{\ge}{2}, q_{0,1},...;\ge) }
     {\tau (q_{0,0}+x-\frac{\ge}{2},q_{0,1},...;\ge)} \\
&&
P_R:=
\tilde w_0 + \frac{\tilde w_1}{\gl}+ \frac{\tilde w_2}{\gl^{2}} +\ldots :=
\frac{\exp\(\sum_{n\geq 0} n!\lambda^{-n-1} \ge\d_{n,1}\)
\tau (q_{0,0}+x+\frac{\ge}{2},q_{0,1},...;\ge) }
     {\tau (q_{0,0}+x-\frac{\ge}{2},q_{0,1},...;\ge)}
\een
We call $\tau$ a {\em tau-function of the ETH} if $\P_L$ and $\P_R$ are
wave operators of the ETH. For a given pair of wave operators the 
tau-function is unique up to multiplication by a non-vanishing function 
independent of $q_{0,0}$ and $q_{n,1}$ with all $n\geq 0$.
In \cite{CDZ} the tau-functions of the ETH were introduced via the 
Lax operator $\L$. The two definitions should agree, up to a factor 
depending on $q_{n,0},\ n>1$ and up to a shift by $\ge/2$ of the 
translation variable $x$.  However, we could not prove this fact.  

Introduce the vertex operators
\begin{align} \notag 
\Gamma^{\pm \ga} & := & \exp \left\{ \pm \frac{1}{2}\sum_{n\geq 0}
\left[ \frac{\gl^{n+1}}{(n+1)!}\frac{q_{n,1}}{\ge}  + 
\frac{2\gl^n}{n!}\( \log {\gl} - \CC_n\) \frac{q_{n,0}}{\ge}
 \right]  \right\} \\ \label{Gamma-alpha} & & 
\times \exp \left\{\mp \frac{\ge}{2}\d_{0,0}  \mp 
\sum_{n\geq 0} \frac{n!}{\gl^{1+n}} \ge \d_{n,1} \right\} .\end{align} 
We remark that the involved functions of $\gl$ are 
derivatives and anti-derivatives of 
\footnote{by an anti-derivative of a function $f$
we mean one of its primitives, i.e. a function $F$ such that $F'=f$} 
$1/2$ and $\gl^{-1}$.

Because of the logarithmic terms, the vertex operators 
$\Gamma^{\pm \ga} \otimes \Gamma^{\mp\ga}$ under the analytic 
continuation around $\gl=\infty$ are multiplied 
by the monodromy factors 
\beq \label{monodromy}
\exp \left\{ \pm \frac{2\pi i}{\ge} \sum_{n\geq 0} \frac{\gl^n}{n!} 
( q_{n,0} \otimes 1 - 1\otimes q_{n,0}) \right\} .
\eeq
As a result, the expressions similar to (\ref{KdV}) do not expand into
Laurent series near $\gl=\infty$. 
To offset the complication we need to generalize the concept of 
vertex operators. 

We will allow vertex operators with coefficient in the algebra 
$\A$ of differential operators $\sum a_i(x,\ge)\d_x^i,$ where each
$a_i(x,\ge)$ is a formal Laurent series in $\ge$ with coefficients infinitely 
differentiable functions in $x$. Let $\#$ be the {\em antiinvolution} 
on $\A$ which acts on the generators $\d_x$ and $x$ as follows: 
\beq \label{involution}(\ge\d_x)^\# = -\ge\d_x+\log {Q}, \ \ x^\#=x .\eeq
Introduce the vertex operator 
(note that $n=0$ is excluded from the summation range)
\ben
\Gamma^{\gd} = \exp \( - \sum_{n>0} 
\frac{\lambda^n}{\ge n!}(\ge\d_x-\log \sqrt{Q}) q_{n,0}\)
                 \ \exp \( x\d_{0,0} \).
\een
We will say that $\tau$ satisfies the {\em Hirota Quadratic Equation
of the ETH} (shortly --- satisfies the HQE) if
\beq \label{HQE}  \frac{d\gl}{\gl} 
\  (\Gamma^{\gd \#}\otimes \Gamma^{\gd})\  
\left( \Gamma^{\ga}\otimes \Gamma^{-\ga}
-\Gamma^{-\ga}\otimes\Gamma^{\ga} \right)\ (\tau \otimes \tau ) \eeq  
\indent {\em computed at $q'_{0,0}-q''_{0,0}=m\ge$
is regular in $\gl$ for each $m\in \Z$.}   

\noindent 
The expression (\ref{HQE}) is interpreted as taking values in the algebra
$\A$ of differential operators with coefficients depending on 
$\q',\q'',\ge$ and $\gl$.  Note that 
\beq\label{monodromy_killer}
\Gamma^{\gd \#}\tensor \Gamma^{\gd} = 
\exp(x\d_{0,0}')\exp\(
\sum_{n>0}\frac{\gl^n}{\ge n!}(\ge\d_x-\log\sqrt Q)(q_{n,0}'-q_{n,0}'') \)
\exp(x\d_{0,0}'').
\eeq
Denote the monodromy factor \eqref{monodromy} by $M$. Moving $M$ from right
to left through the three exponential factors in \eqref{monodromy_killer}
we get
\ben
&&
\(\Gamma^{\gd \#}\otimes \Gamma^\gd \) M = 
M \exp(x\d_{0,0}') 
\exp\(\pm\frac{2\pi i}{\ge}x\)\\ 
&&
\exp\(
\sum_{n>0}\frac{\gl^n}{\ge n!}(\ge\d_x-\log\sqrt Q)(q_{n,0}'-q_{n,0}'') \)
\exp\( \mp\frac{2\pi i}{\ge}x\)\exp(x\d_{0,0}'')=\\
&&
M\exp(x\d_{0,0}') \exp\(
\sum_{n>0}\frac{\gl^n}{\ge n!}
(\ge(\d_x-\(\pm2\pi i/\ge\))-\log\sqrt Q)(q_{n,0}'-q_{n,0}'') \)
\exp(x\d_{0,0}'')= \\
&&
M\exp\(\mp\frac{2\pi i}{\ge}
\sum_{n>0}\frac{\gl^n}{ n!}(q_{n,0}'-q_{n,0}'') \)
\(\Gamma^{\gd \#}\otimes \Gamma^\gd \)=
e^{\pm \frac{2\pi i}{\epsilon}(q'_{0,0}-q''_{0,0})}\ 
\(\Gamma^{\gd \#}\otimes \Gamma^\gd \). 
\een
Thus when $q'_{0,0}-q''_{0,0} \in \Z\ge $, the expression 
(\ref{HQE}) is single-valued near $\gl=\infty$. After the change 
$\y=(\q'-\q'')/2, \x=(\q'+\q'')/2$ it expands (for each $m$) 
as a power series in $\y$ ($y_{0,0}=m\ge $ excluded) with coefficients 
which are Laurent series in $\gl^{-1}$ 
(whose coefficients are {\em differential operators in $x$}
depending on $\x$ via $\tau$, its translations and partial derivatives). 

\begin{theorem}\label{t1}
A non-vanishing function $\tau$  is a tau-function of the Extended Toda 
Hierarchy if and only if it satisfies the Hirota Quadratic Equations
\eqref{HQE}.
\end{theorem}

As a by-product we obtain a Hirota equation for the (unextended) 
Toda Lattice Hierarchy \cite{UT} described in the Lax form by (\ref{laxeq1}).
The concept of tau-functions
easily carries over to this case. Introduce the vertex operators
\ben \Gamma^{\pm \gb} = \exp \left\{ \pm \frac{1}{2} \sum_{n\geq 0}
\frac{\gl^{n+1}}{(n+1)!}\frac{q_{n,1}}{\ge} \right\} \ 
\exp \left\{\mp \frac{\ge}{2}\d_{0,0}  \mp 
\sum_{n\geq 0} \frac{n!}{\gl^{1+n}} \ge \d_{n,1} \right\} .\een 

\begin{corollary}\label{c1} A non-vanishing function 
$\tau$ of $(q_{0,0};q_{1,0},q_{1,1},...;\ge)$ is a tau-function of the 
Toda Lattice Hierarchy (\ref{laxeq1}) 
if and only if for each $m\in \Z$
\ben \frac{d\gl}{\gl} \left\{ \(\frac{\gl}{\sqrt{Q}}\)^m 
\Gamma^{\gb}\otimes \Gamma^{-\gb} - \(\frac{\gl}{\sqrt{Q}}\)^{-m} 
\Gamma^{-\gb}\otimes \Gamma^{\gb} \right\} (\tau \otimes \tau) \een
\indent computed at $q_{0,0}'-q_{0,0}''=m\ge$ is regular in $\gl$.
\end{corollary} 

{\bf Acknowledgments.} I am thankful to A. Givental for the many stimulating 
discussions. I am thankful also to the referees for their valueable remarks 
which helped me to improve the exposition and to fix some inaccuracies in 
the paper.

\sectionnew{Wave operators and tau-functions of the ETH}

The arguments in this section are parallel to the ones in \cite{UT}.
\subsection{Wave operators} For any Lax operator $\L$ put
\ben
\A_{n,0}=\( \frac{2\L^n}{n!\epsilon}(\log \L -\CC_n)  \)_+\ ,\
&&
\A_{n,1}=\(\frac{\L^{n+1}}{(n+1)!\epsilon}\)_+ \ .
\een  
The compatibility of \eqref{bn1}--\eqref{an0} is equivalent to 
the Zakharov-Shabat equations:
\begin{lemma}\label{z-sh}
If $\L$ satisfies the equations of the ETH then 
\ben
\d_{n,\ga}\A_{m,\gb}-\d_{m,\gb}\A_{n,\ga} = [\A_{n,\ga},\A_{m,\gb}]\ , 
\een
where $m,n\geq 0$ and $\ga,\gb=0,1$. 
\end{lemma}
\proof
The proof that a solution to the Lax equations implies the Zakharov-Shabat
equations is standard: see for example \cite{PvM} Theorem 1.1, where 
this is done for the KP hierarchy. However, we remark that in 
our case one should use the following non-trivial formula, 
which is proven in \cite{CDZ}:
\ben
\d_{n,\ga} \log \L = [\A_{n,\ga},\log \L], \ \ga=0,1.
\een
\qed

As a corollary we find that the Cauchy problem for
the system of differential equations \eqref{bn1}--\eqref{an0} has 
a unique solution.
\begin{proposition}\label{ETH-w}
A Lax operator $\L$ satisfies the equations of the  ETH if and only if
there is a pair of dressing operators $\P_L$ and $\P_R$, which satisfies
\eqref{bn1}--\eqref{an0}.
\end{proposition}
\proof
Fix some sequence of times $\q^0$. Consider the Lax operator $\L^0:=\L(\q^0)$.
There exist dressing operators $\P_L^0$ and $\P_R^0$, i.e.
\ben
\L^0 = \P_L^0\Lambda(\P_L^0)^{-1}=((\P_R^0)^{-1}\Lambda\P_R^0)^\#.
\een
Let $\P_L$ and $\P_R$ be solutions to the system of equations
\eqref{bn1}--\eqref{an0}, satisfying the initial conditions 
$\P_L|_{\q=\q^0}=\P_L^0$ and $\P_R|_{\q=\q^0}=\P_R^0$.  One checks 
that the operators
\ben
\L\P_L-\P_L\Lambda \mbox{ and  } 
\L\P_R^\# - Q\P_R^\#\Lambda^{-1}  
\een
also satisfy \eqref{bn1}--\eqref{an0}. Since at $\q=\q^0$ both operators are
0 we obtain that they are identically zero. This proves that $\P_L$ and 
$\P_R$ are wave operators. 

In the other direction we need to use only that if $\L=\P_L\Lambda\P_L^{-1}$
then
\ben
\d_{n,\ga}\L = [(\d_{n,\ga}\P_L)\P_L^{-1},\L],\ n\geq 0,\ \ga=0,1.
\een
\qed

\subsection{Characterization of the wave operators of the ETH} 
Let $\A[[\q]]$ be the algebra of formal power series in $\q$ with 
coefficients in the algebra $\A$ of differential operators. We identify
the Lax and the corresponding dressing operators with vectors in the
space $\A[[\q]][[\Lambda^{\pm1}]]$ of formal series in $\Lambda$, where
the translation operator $\Lambda$ is identified with a formal symbol
which commutes with the elements of $\A[[\q]]$ in the same ways as
$e^{\ge\d_x}$ does, i.e. if 
$\sum_{k\geq 0} a_k(x,\q;\ge)\d_x^k \in \A[[\q]]$ then
$$
\Lambda^{\pm 1}\,\(\sum_{k\geq 0} a_k(x,\q;\ge)\d_x^k\) = 
         \(\sum_{k\geq 0} a_k(x\pm\ge,\q;\ge)\d_x^k\)\, \Lambda^{\pm 1}.
$$
 
Let
\ben
\Q=\sum_{k\in \Z} b_k\Lambda^k = \sum_{k\in \Z}\Lambda^k\tilde b_k
\in \A[[\q]][[\Lambda^{\pm1}]].
\een
Then the series 
\ben
\sum_{k\in \Z} b_k\lambda^k
\mbox{ and }
\sum_{k\in \Z}\tilde b_k\gl^k
\een
will be called respectively {\em left} and {\em right symbols} of $\Q$.

Let $\P_L$ and $\P_R$ be two arbitrary operator series of the form
respectively \eqref{dressP} and \eqref{dressQ}. We will assume that the 
coefficients of $\P_L$ and $\P_R$ are formal series in 
$x+q_{0,0},\ q_{0,1},\ q_{n,i},\ n>0 ,\ i=0,1.$ Note that the wave 
operators of the ETH also have this form. Our goal is to see what further 
restrictions should be imposed on $\P_L$ and $\P_R$ so that they become 
wave operators of the ETH. Introduce the following two series
in $\A[[\q]][[\Lambda^{\pm1}]]$:
\beqa
\label{wl}
&&
\W_L(x,\q,\Lambda) =
\P_L \exp\left\{ 
\sum_{n\geq 0}\frac{\Lambda^{n+1}}{2\ge(n+1)!}q_{n,1}
+\sum_{n> 0}\frac{\Lambda^n}{\ge n!}
(\ge \d_x-\CC_n)q_{n,0}   \right\},    \\
\label{wr}
&&
\W_R(x,\q,\Lambda) = 
 \exp\left\{ 
-\sum_{n\geq 0}\frac{ \Lambda^{n+1} }{2\ge(n+1)!}q_{n,1}
-\sum_{n> 0}\frac{ \Lambda^{n}   }{\ge n!}
(\ge \d_x-\CC_n)q_{n,0}   \right\}\P_R.
\eeqa
Finally, denote by $W_L$ the left symbol of $\W_L$ and by
$W_R$ the right symbol of $\W_R$.  
\begin{proposition}\label{wave-operators}
Let $\q'$ and $\q''$ be such that $q_{0,0}'=q_{0,0}''$. 
The following conditions are equivalent:
\item{ (a)} 
$\P_L$ and $\P_R$ are wave operators of the ETH .
\item{ (b)} 
$\P_L$ and $\P_R$ satisfy the following identities:
\beqa
&&
\label{HQE-o1}
\P_L\Lambda\P_L^{-1} = (\P_R^{-1} \Lambda \P_R)^\# \\
&&
\label{HQE-o0}
\W_L(x,\q',\Lambda)\W_R(x,\q'',\Lambda) = 
\{ \W_L(x,\q'',\Lambda)\W_R(x,\q',\Lambda)\}^\# .
\eeqa
\item{ (c)} 
For all integers $r\geq 0$ the following identity holds:
\beqa
\label{HQE-o}
\W_L (x,\q',\Lambda) \Lambda^r \W_R(x,\q'',\Lambda) =
\left\{
\W_L(x,\q'',\Lambda) \Lambda^r \W_R (x,\q',\Lambda)
\right\}^\# .
\eeqa
\item{ (d)}
For all integers $m$ and $r\geq 0$ the following identity holds
\footnote{The residue here is interpreted  
as the coefficient in front of $\gl^{-1}.$}: 
\beqa
&&
\label{HQE-s}
\res_{\lambda = \infty} \left\{  
\lambda^r \( \frac{\gl}{\sqrt Q}\)^m\ 
W_L(x,\q',\gl)   W_R(x-m\ge,\q'', \gl)  \right\} 
\frac{d\gl}{\gl}  = \\ \notag
&&
\res_{\gl = \infty} 
\left\{ \gl^r\( \frac{\gl}{\sqrt Q} \)^{-m}
\(W_R( x,\q',\gl )\)^\#  \(W_L(x-m\ge,\q'',\gl)\)^\# \right\} 
\frac{d\gl}{\gl} .
\eeqa
\end{proposition} 
\proof

\boxed{(d) \Leftrightarrow (c)}
Let $m$ and $r\geq 0$ be some arbitrary integers and $q_{0,0}'=q_{0,0}''$. 
Put 
\ben
\W_L(x,\q,\Lambda)=\sum_{i\in \Z} a_i(x,\q, \d_x)\Lambda^i 
\mbox{ and }
\W_R(x,\q,\Lambda)=\sum_{j\in \Z} \Lambda^jb_j(x,\q, \d_x)
\een
and compare the coefficients in front of $\Lambda^{-m}$ in 
\eqref{HQE-o}: 
\ben
\sum_{i+j=-m-r} a_i(x,\q', \d_x)b_j(x-m\ge,\q'', \d_x) = 
\sum_{i+j=m-r} Q^{m} b_j^\#(x,\q', \d_x)a_i^\#(x-m\ge,\q'', \d_x). 
\een
This equality can be written also as
\ben
&&
\res_{\lambda = \infty} \left\{  
\lambda^r \( \frac{\gl}{\sqrt Q}\)^m\ 
W_L(x,\q',\gl)   W_R(x-m\ge,\q'', \gl)  \right\} 
\frac{d\gl}{\gl}  = \\ \notag
&&
\res_{\gl = \infty} 
\left\{ \( \gl/Q \)^{-r}\(\frac{\gl}{\sqrt Q}\)^{m}
\(W_R( x,\q',Q\gl^{-1} )\)^\#  \(W_L(x-m\ge,\q'',Q\gl^{-1})\)^\# \right\} 
\frac{d\gl}{\gl}. 
\een
The R.H.S. of the last equality is the coefficient in front
of $\gl^0$ in the series inside the $\{\ \}$-brackets. Replacing 
$\gl$ with $Q/\gl$ inside the $\{\ \}$-brackets we get  \eqref{HQE-s}.

\boxed{(c)\Rightarrow (b)}
Assume first that \eqref{HQE-o} holds. The identity \eqref{HQE-o0} is 
obtained by letting $r=0$. To prove \eqref{HQE-o1},  
put $r=1$ and $\q'=\q''$ in \eqref{HQE-o}. We get the equality 
$\P_L\Lambda\P_R = (\P_L\Lambda\P_R)^\#$ which can be written also as:
\ben
\P_L\Lambda\P_L^{-1} (\P_L\P_R) = (\P_R^{-1} \Lambda \P_R)^\# (\P_L\P_R)^\#
\een
It remains only to notice that $\P_L\P_R=(\P_L\P_R)^\#.$ Indeed, 
in \eqref{HQE-o} let $r=0$ and $\q'=\q''$.

\boxed{(b)\Rightarrow (c)}
Put 
$\L=\P_L\Lambda\P_L^{-1} = (\P_R^{-1} \Lambda \P_R)^\#$. Note that 
$\L$ should have the form $\Lambda+ a_0+a_1\Lambda^{-1}$. Thus $\L$ is
a Laurent polynomial in $\Lambda$  and it makes sense to multiply 
by $\L$ any $\Lambda$-series (possibly infinite in both directions). 
In particular multiply  
\eqref{HQE-o0} by $\L^r,\ r\geq 0$, and use 
\ben
\L \W_L = \W_L\Lambda,
\quad
\W_R \L^\# = \Lambda\W_R 
\een
to obtain \eqref{HQE-o}.

\boxed{(b) \Rightarrow (a)}
Note that \eqref{HQE-o1}
implies that $\L:=\P_L\Lambda\P_L^{-1}$ is a Lax operator.  
Let's prove \eqref{an0} and \eqref{bn0}. When $n=0$ we need to show 
that  $\d_{0,0}\P_L=\d_x\P_L$ and $\d_{0,0}\P_R^\#=\d_x\P_R^\# $. 
Both are satisfied by definition.
Assume that $n>0$. Differentiate \eqref{HQE-o0} with respect to
$q_{n,0}'$ and  then put $\q'=\q''$:
\ben
&&
(\d_{n,0}\P_L)\P_L^{-1}+\P_L\frac{\Lambda^n}{\ge n!}
(\ge\d_x-\log \sqrt Q - \CC_n)\P_L^{-1} = \\
&&
\(\P_R^{-1}\d_{n,0}\P_R \)^\# 
-\( \P_R^{-1}\frac{\Lambda^n}{\ge n!}
(\ge\d_x-\log \sqrt Q - \CC_n)\P_R \)^\#
\een
Using the definition of $\log \L$ and \eqref{HQE-o1} the last identity
simplifies to 
\ben
(\d_{n,0}\P_L)\P_L^{-1} +\frac{2\L^n}{\ge n!}
( \log \L -\CC_n) = \(\P_R^{-1}\d_{n,0}\P_R\)^\#
\een
Since $(\d_{n,0}\P_L)\P_L^{-1}$ contains only negative powers of 
$\Lambda$ and $\(\P_R^{-1}\d_{n,0}\P_R\)^\#$ -- non-negative,
we get \eqref{an0} and \eqref{bn0} by separating the negative and the 
positive part of the equation. 

To check that \eqref{an1} and \eqref{bn1} hold:
differentiate \eqref{HQE-o} with respect to $q_{n,1}'$, 
put $\q'=\q''$ and then apply a similar argument.  Thus $\P_L, \ \P_R$ 
is a pair of wave operators.

\boxed{(a)\Rightarrow (b)}
Let $\ga=(\ga_{0,1},\ga_{1,0},\ga_{1,1},\ga_{2,0},\ga_{2,1},\ldots )$ 
be a multindex with only finetly many non-zero components.
Put 
$$
\d^\ga:=\d_{0,1}^{\ga_{0,1}}\d_{1,0}^{\ga_{1,0}}\d_{1,1}^{\ga_{1,1}}
\d_{2,0}^{\ga_{2,0}}\d_{2,1}^{\ga_{2,1}}\ldots\ ,
$$
where $\d_{n,i}=\d/\d q_{n,i}$ (note that the differentiation 
$\d/\d q_{0,0}$ is not involved). First, using induction on $|\ga |$,  
we prove that 
\beq
\label{HQE-oa}
\(\d^\ga\W_L(x,\q,\Lambda)\)\W_R(x,\q,\Lambda) = 
\{\W_L(x,\q,\Lambda)\d^\ga \W_R(x,\q,\Lambda)\}^\#,
\eeq
When $\ga=0$ we need to prove that 
\ben
\W_L(x,\q,\Lambda) \W_R(x,\q,\Lambda) = 
\(\W_L(x,\q,\Lambda) \W_R(x,\q,\Lambda)\)^\#,
\een
which is equivalent to $\P_L\P_R=(\P_L\P_R)^\#$. We claim that the last 
identity holds for any pair $(\P_L,\P_R)$  of dressing 
operators corresponding to the Lax operator 
$\L=\Lambda +u +Qe^v\Lambda^{-1}$. 
Note that  $\L^\# =\tilde w_0^{-1}\L \tilde w_0$, where 
$\tilde w_0$  is such that $e^{v(x)}=\tilde w_0(x)/\tilde w_0(x-\ge)$. Thus
$(\P_L,\widetilde \P_R:=\P_L^{-1}\tilde w_0)$ is another pair of
dressing operators for $\L$. By the uniqueness of the dressing operators
we get that there is a series $K=a_0+a_1\Lambda^{-1}+\ldots$, with 
coefficients $a_i$ independent of $x$ such that $\P_R=K\tilde \P_R$.  Thus
\ben
&&
\P_L\P_R = \P_L(a_0+a_1\Lambda^{-1} +\ldots )\P_L^{-1}\tilde w_0 = 
(a_0+a_1\L^{-1} +\ldots )\tilde w_0 = \\
&&
\tilde w_0 (a_0+a_1(\L^\#)^{-1}+\ldots)=
\( (a_0+a_1\L^{-1} +\ldots )\tilde w_0 \)^\# = (\P_L\P_R)^\#. 
\een
Assume now that  \eqref{HQE-oa} is true for some $\ga$ and 
differentiate in $\d_{n,0}$. We need only to check that
\beq
\label{induction-step}
\d^\ga\W_L\d_{n,0}\W_R = \{\d_{n,0}\W_L\d^\ga\W_R\}^\#.
\eeq
Note that
\ben
&&
\d_{n,0}\W_L = \( (\d_{n,0}\P_L)\P_L^{-1}+
\P_L\frac{\Lambda^n}{\ge n!}( \ge\d_x-\log \sqrt Q-\CC_n )\P_L^{-1}\)\W_L, \\
&&
\d_{n,0}\W_R = \W_R\( \P_R^{-1}\d_{n,0}\P_R-
\P_R^{-1}\frac{\Lambda^n}{\ge n!}( \ge\d_x-\log \sqrt Q-\CC_n )\P_R\).
\een
Now \eqref{induction-step} follows from \eqref{bn0}, \eqref{an0} and the 
inductive assumption \eqref{HQE-oa}. Thus if we increase $\ga_{n,0}$ by 1
then \eqref{HQE-oa} still holds. Similarly if we increase $\ga_{n,1}$ by 1 
then \eqref{HQE-oa} still holds. The induction is completed. 

Using the Taylor's formula expand both sides of \eqref{HQE-o0} 
about $\q'=\q''$. Then the coefficients in front of $(\q'-\q'')^\ga$ 
are equal exactly when \eqref{HQE-oa} holds. 
\qed

\subsection{Tau-functions}
Let 
\ben 
&&
\P_L(x,\q,\Lambda)=
1+w_1(x,\q;\ge)\Lambda^{-1}+w_2(x,\q;\ge)\Lambda^{-2}+\dots ,\\
&&
\P_R(x,\q,\Lambda)=
\tilde w_0(x,\q;\ge) + \Lambda^{-1}\tilde w_1(x,\q;\ge) + \Lambda^{-2} 
\tilde w_2(x,\q;\ge) +\ldots 
\een
be a pair of wave operators of the ETH. Denote by:
\ben
&&
P_L(x,\q,\lambda):=
1+w_1(x,\q;\ge)\gl^{-1}+w_2(x,\q;\ge)\gl^{-2}+\ldots \\
&&
P_R(x,\q,\lambda):=
\tilde w_0(x,\q;\ge) + \tilde w_1(x,\q;\ge)\gl^{-1}+ 
                   \tilde w_2(x,\q;\ge)\gl^{-2} +\ldots
\een
the left and the right symbol respectively of $\P_L$ and $\P_R.$
For shortness denote by $[\gl^{-1}]$ the sequence 
of times with components
\ben
[\gl^{-1}]_{n,0}=0,\ 
[\gl^{-1}]_{n,1}=n!\gl^{-n-1}\ge.
\een

\begin{lemma}{\label{identities}}
The following identities hold
\beqa
&&
\label{id1}
P_L(x,\q,\gl)P_R(x-\ge,\q-[\gl^{-1}],\gl) = \tilde w_0(x-\ge,\q-[\gl^{-1}]),\\
&&
\label{id2}
P_L(x,\q,\gl_1)P_R(x-\ge,\q-[\gl_1^{-1}]-[\gl_2^{-1}],\gl_1)= \\ \notag
&&
\ \ \ \ \ \ 
=P_L(x,\q,\gl_2)P_R(x-\ge,\q-[\gl_1^{-1}]-[\gl_2^{-1}],\gl_2),  \\
&&
\label{id4}
P_L(x,\q,\gl)P_R(x,\q-[\gl^{-1}],\gl) = \tilde w_0 (x,\q).
\eeqa
\end{lemma}
\proof
Let $\q'$ and $\q''$ be two sequences of time variables such that 
$q'_{n,0}=q''_{n,o},\ n\geq 0$. 
The identities \eqref{id1}--\eqref{id4} are consequence of the 
following one: 
\beqa
\label{specialize2}
&& 
\P_L(x,\q',\Lambda)\exp\Big(\sum_{n\geq 0}\frac{\Lambda^{n+1}}{\ge (n+1)!}
(q_{n,1}'-q_{n,1}'')\Big) \P_R(x,\q'',\Lambda) = \\\notag
&&
\P_R^\#(x,\q',\Lambda)\P_L^\#(x,\q'',\Lambda).
\eeqa
The proof of \eqref{specialize2} is completely analogous to the argument 
in the implication 
$a)\Rightarrow b)$ in \prref{wave-operators} and it will be omitted. 

To prove \eqref{id1}: in \eqref{specialize2} put 
$\q''=\q'-[\gl^{-1}]$. The exponential factor turns into
\ben
\exp \Big( \sum_{n\geq 0} \frac{(\lambda^{-1}\Lambda)^{n+1}}{n+1}
\Big) = 
 (1-\lambda^{-1}\Lambda)^{-1}=\sum_{N\geq 0}(\gl^{-1}\Lambda)^N.
\een 
Comparing the coefficients in front of $\Lambda^{-1}$ we find
\ben
&&
\sum_{N\geq 0, j+k=N+1} w_j(x,\q')\gl^{-N}\Lambda^{-1}\tilde w_k(x,\q'')=0, \\
&&
\sum_{N\geq 0, j+k=N+1} w_j(x,\q')\gl^{-j}\gl^{-k}\tilde w_k(x-\ge,\q'')=0,\\
&&
P_L(x,\q',\lambda)P_R(x-\ge,\q'-[\gl^{-1}],\gl)
-\tilde w_0(x-\ge,\q'-[\gl^{-1}])=0.
\een

To prove \eqref{id2}: in  \eqref{specialize2} put 
$\q''=\q'-[\gl_1^{-1}]-[\gl_2^{-1}]$. The exponential factor turns into
\ben
(1-\gl_1^{-1}\Lambda)^{-1}(1-\gl_2^{-1}\Lambda)^{-1} = 
\frac{\gl_1\gl_2}{\gl_2-\gl_1}\{(1-\gl_1^{-1}\Lambda)^{-1}-
(1-\gl_2^{-1}\Lambda)^{-1}\}\Lambda^{-1}.
\een 
Comparing the coefficients in front of $\Lambda^{-1}$ we arrive at
\eqref{id2}.

To prove \eqref{id4}: in \eqref{specialize2} put $\q''=\q'-[\gl^{-1}]$ and 
compare the coefficients in front of $\Lambda^0$. 
\qed

The main result in this subsection is the following proposition:
\begin{proposition}\label{tau-function}
Given a pair of wave operators $\P_L$ and $\P_R$ of the ETH there exists 
a corresponding tau-function, which is unique up to multiplication by a 
non-vanishing function independent of $q_{0,0}$ and $q_{n,1},\ n\geq 0$.
\end{proposition}
\proof
We need to prove that the system of equations
\footnote{By definition, the tau-functions corresponding to the 
wave operators $\P_L$ and $\P_R$ are $\tau(x-\ge/2,\q)$, where 
$\tau(x,\q)$ is a solution to the system \eqref{t_eq_1}--\eqref{t_eq_3}.}
\beqa\label{t_eq_1}
&&
P_L(x,\q,\gl)=\frac{\tau(x,\q-[\gl^{-1}])}{\tau(x,\q )}, \\ \label{t_eq_2} 
&&
P_R(x,\q,\gl)=\frac{\tau(x+\ge,\q+[\gl^{-1}])}{\tau(x,\q )}, \\ \label{t_eq_3}
&&
\d_{0,0}  \tau(x,\q) = \d_x \tau(x,\q).
\eeqa
has a solution $\tau(x,\q)$ unique up to  multiplication by a 
non-vanishing function independent of $q_{0,0}$ and $q_{n,1},\ n\geq 0$.
The system is equivalent to: 
\beqa
&&
\label{eq1}
\log P_L = \(e^{-\sum n!\gl^{-n-1}\ge\d_{n,1}}-1\)\log \tau, \\
\label{eq2}
&&
\log P_R = \(e^{\ge\d_{x}+\sum n!\gl^{-n-1}\ge\d_{n,1}}-1\)\log \tau, \\
\label{eq3}
&&
\d_{0,0} \log \tau(x,\q) = \d_x \log \tau(x,\q).
\eeqa
Expanding in the powers of $\gl$ we find:
\ben
&&
\log P_L = \sum_{N\geq 1} b_N(x,\q)\gl^{-N}, \\
&&
\log P_R = \sum_{N\geq 0} \tilde b_N(x,\q)\gl^{-N}
\een
and 
\beq
\label{schur}
\exp\({-\sum n!\gl^{-n-1}\ge\d_{n,1}}\)-1 = 
\sum_{N\geq 1}a_N(\d_{0,1}, \d_{1,1},\ldots )\gl^{-N}.
\eeq
Comparing the coefficients in front of the powers of $\gl$ we 
get the following system of partial differential equations (note that
\eqref{eq2} is equivalent to \eqref{tau-de1} and \eqref{tau-de3}):
\beqa 
\label{tau-de1}
&&
 a_N(-\d_{0,1}, -\d_{1,1},\ldots ) \log\tau =  
e^{-\ge\d_x}\tilde b_N(x,\q), \  N\geq 1 ,\\
\label{tau-de2}
&&
a_N(\d_{0,1}, \d_{1,1},\ldots )\log \tau = b_N(x,\q),\ N\geq 1 ,\\
\label{tau-de3}
&&
\ge\d_x\log \tau = \frac{\ge\d_x}{e^{\ge\d_x}-1}\tilde b_0,\\
\label{tau-de4}
&&
\d_{0,0}\log \tau = \d_x \log \tau,
\eeqa
where 
\ben
\frac{\ge\d_x}{e^{\ge\d_x}-1}=\sum_{k=0}^\infty B_k \frac{(\ge\d_x)^k}{k!}, 
\ \ B_k \mbox{ are the Bernoulli numbers}.
\een
Let us exclude the first equation. 
Later on we will see that it is a consequence from the rest.  
Note that the differential operators $a_1(\d),a_2(\d),\ldots $ 
generate polynomially the ring $\C[\d_{0,1}, \d_{1,1},\ldots ]$. Thus the 
system \eqref{tau-de2}--\eqref{tau-de4} can be written in the form 
\beqa
\label{tau-de}
\begin{cases}
\vspace{.1in}
\d_{n,1}\log \tau = \beta_n(x,\q), \ n\geq 1, \\
\vspace{.1in}
\ge\d_x\log \tau = \frac{\ge\d_x}{e^{\ge\d_x}-1}\tilde b_0,\\
\d_{0,0}\log \tau = \d_x \log \tau,
\end{cases}
\eeqa
where $\beta_n(x,\q)$ depend polynomially  on $b_n(x,\q)$  and their 
derivatives.  
We need to show that the system is compatible. 

Since $P_L$ and $P_R$ depend on $x+q_{0,0}$ the equation on the third 
line of \eqref{tau-de} is compatible with the rest of the equations.

The  equation on the second line of \eqref{tau-de} is compatible with the 
equations on the first line if and only if
\ben
a_i(\d)\frac{\ge\d_x}{e^{\ge\d_x}-1}\tilde b_0 =\ge\d_x b_i,\ i\geq 1,
\een
which is  equivalent to
\ben
a_i(\d)\tilde b_0=(e^{\ge\d_x}-1)b_i,\ i\geq 1.
\een
Let us write a generating series for these identities
\ben
&&
\sum_{i\geq 1}  a_i(\d)\tilde b_0 \gl^{-i} =
\sum_{i\geq 1} (e^{\ge\d_x}-1)b_i \gl^{-i}. 
\een
Comparing with the expansions of $\log P_L$, \eqref{schur} and letting
$\tilde b_0=\log \tilde w_0$ we get
\ben
\log \frac{\tilde w_0(x,\q-[\gl^{-1}])}{\tilde w_0(x,\q)}=
\log \frac{P_L(x+\ge,\q,\gl)}{P_L(x,\q,\gl)},
\een
which is equivalent to
\beq
\label{identity-b}
\tilde w_0(x,\q-[\gl^{-1}])P_L(x,\q,\gl)=
\tilde w_0(x,\q) P_L(x+\ge,\q,\gl).
\eeq
Similarly, the compatibility between the equations on
the first line  of \eqref{tau-de} is equivalent  to
\beqa
\label{identity-a}
P_L(x,\q,\gl_1)P_L(x,\q-[\gl_1^{-1}],\gl_2)=
P_L(x,\q,\gl_2)P_L(x,\q-[\gl_2^{-1}],\gl_1). 
\eeqa
Thus the compatibility of \eqref{tau-de} is equivalent to \eqref{identity-b}
and \eqref{identity-a}. 

We will show that  \eqref{identity-b}
and \eqref{identity-a} follow from the identities in \leref{identities}.
To prove \eqref{identity-b}: divide \eqref{id2} by 
$\tilde w_0 (x-\ge,\q-[\gl_1^{-1}]-[\gl_2^{-1}])$ and then apply
\eqref{id1}. To prove \eqref{identity-a}:  substitute $x$ 
with $x-\ge$ in \eqref{id1} and use \eqref{id4}.  

Thus the system \eqref{tau-de} is compatible. Let $\tau$ be 
a solution -- it is unique up to multiplication
by a constant independent of $x, q_{0,0}$ and $q_{n,1},\ n\geq 0$. 
Then we have:
\ben
P_R(x,\q,\gl) = \frac{\tau(x+\ge,\q+[\gl^{-1}])}{\tau(x,\q)},\ \ 
\tilde w_0 = \frac{\tau(x+\ge,\q)}{\tau(x,\q)}.
\een
Using \eqref{id4} we get
\ben
P_L(x,\q,\gl) = \frac{\tau(x,\q-[\gl^{-1}])}{\tau(x,\q)}.
\een
\qed

\subsection{Proof of \thref{t1}}
Assume that $\tau$ is a non-vanishing function and let $\P_L$ and $\P_R$ 
be the
corresponding operators (see the introduction). It is enough to 
prove that the HQE are equivalent to condition (d) in \prref{wave-operators}. 

After a straightforward computation one finds
{\allowdisplaybreaks 
\begin{align} \label{f1}
\Gamma^{\gd \#} \Gamma^{\ \ga}\tau & = 
\tau(x-\ge/2,\q)  \(\frac{\gl}{\sqrt{Q}} \)^{\ q_{0,0}/\ge}
W_L(x,\q,\gl )\(\frac{\gl}{\sqrt{Q}} \)^{      x/\ge} \\ \label{f2}
\Gamma^{\gd \#} \Gamma^{-\ga}\tau  & = 
\tau(x-\ge/2,\q)  \(\frac{\gl}{\sqrt{Q}} \)^{-q_{0,0}/\ge}
\(W_R(x,\q,\gl )\)^\#\(\frac{\gl}{\sqrt{Q}} \)^{      -x/\ge} \\ \label{f3}
\Gamma^{\gd} \Gamma^{-\ga}\tau  & = 
\(\frac{\gl}{\sqrt{Q}} \)^{-q_{0,0}/\ge}
\(\frac{\gl}{\sqrt{Q}} \)^{      -x/\ge}
W_R(x,\q,\gl)\ 
\tau(x-\ge/2,\q)  \\ \label{f4}
\Gamma^{\gd} \Gamma^{\ \ga}\tau & = 
\(\frac{\gl}{\sqrt{Q}} \)^{\ \ q_{0,0}/\ge}
\(\frac{\gl}{\sqrt{Q}} \)^{  \  \    x/\ge}
\(W_L(x,\q,\gl)\)^\#\ 
\tau(x-\ge/2,\q) .
\end{align}}
Let us prove the first identity. The other three are derived
in a similar way. 
{\allowdisplaybreaks 
\ben
&&
\Gamma^{\gd \#} \Gamma^{\ \ga}\tau =
\exp\(x\d_{0,0}\)
\exp\(\sum_{n>0}\frac{\gl^n}{\ge n!}(\ge\d_x-\log\sqrt Q)q_{n,0}\)\times \\
&&
\times\exp \left\{ \frac{1}{2}\sum_{n\geq 0}
\left[ \frac{\gl^{n+1}}{(n+1)!}\frac{q_{n,1}}{\ge}  + 
\frac{2\gl^n}{n!}\( \log {\gl} - \CC_n\) \frac{q_{n,0}}{\ge}
 \right]  \right\}\\
&&
\times \exp \left\{- \frac{\ge}{2}\d_{0,0}  - 
\sum_{n\geq 0} \frac{n!}{\gl^{1+n}} \ge \d_{n,1} \right\}
\tau(\q;\ge).
\een}
We move the third line of the last formula between the two exponents
on the first line (note that the third line commutes with the preceding two
exponential factors). Then we move the  operator $\exp(x\d_{0,0})$
from left to right. We get
{\allowdisplaybreaks 
\ben
&&
\Gamma^{\gd \#} \Gamma^{\ \ga}\tau =\tau(x-\ge/2,\q)P_L(x,\q,\gl)
\exp\(\sum_{n>0}\frac{\gl^n}{\ge n!}
(\ge(\d_x-\d_{0,0})-\log\sqrt Q)q_{n,0}\) \\
&&
\exp \left\{ \frac{1}{2}\left[\sum_{n\geq 0}
 \frac{\gl^{n+1}}{(n+1)!}\frac{q_{n,1}}{\ge}  + \sum_{n>0}
\frac{2\gl^n}{n!}\( \log {\gl} - \CC_n\) \frac{q_{n,0}}{\ge}\right]
+(\log \gl -\log\sqrt Q) (q_{0,0}+x)/\ge
 \right\}=\\
&&
\tau(x-\ge/2,\q)P_L(x,\q,\gl)
\exp\(\sum_{n>0}\frac{\gl^n}{\ge n!}
(\ge\d_x-\log\sqrt Q)q_{n,0}\) \\
&&
\exp \left\{ \frac{1}{2}\left[\sum_{n\geq 0}
 \frac{\gl^{n+1}}{(n+1)!}\frac{q_{n,1}}{\ge}  + \sum_{n>0}
\frac{2\gl^n}{n!}\( \log {\gl} - \CC_n\) \frac{q_{n,0}}{\ge}\right] +\right.\\
&&
\left.
+(\log \gl-\log\sqrt Q) \[q_{0,0}-\(\sum_{n>0}\frac{\gl^n}{n!}q_{n,0}\)+x\]/\ge
 \right\}=
\tau(x-\ge/2,\q)P_L(x,\q,\gl) \times \\
&& \exp \left\{ \frac{1}{2}\left[\sum_{n\geq 0}
 \frac{\gl^{n+1}}{(n+1)!}\frac{q_{n,1}}{\ge} + 
\sum_{n>0}\frac{\gl^n}{n!}
(\ge\d_x-\log\sqrt Q)\frac{q_{n,0}}{\ge} \right]\right\} 
\(\frac{\gl}{\sqrt Q}\)^{(q_{0,0}+x)/\ge}.
\een}

After substituting formulas \eqref{f1}--\eqref{f4} into 
the HQE we find:
\ben
&&
\left\{ \Gamma^{\gd\#}\Gamma^{\ga}\tau \tensor 
\Gamma^{\gd}\Gamma^{-\ga}\tau -
\Gamma^{\gd\#}\Gamma^{-\ga}\tau \tensor 
\Gamma^{\gd}\Gamma^{\ga}\tau \right\} \frac{d\gl}{\gl} = \\
&&
\left\{ \tau(x-\ge/2,\q') \(\frac{\gl}{\sqrt Q}\)^{(q_{0,0}'-q_{0,0}'')/\ge} 
W_L(x,\q',\gl)
W_R(x,\q'',\gl)\ \tau(x-\ge/2,\q'')  - \right. \\
&&
\left. 
\tau(x-\ge/2,\q') \(\frac{\gl}{\sqrt Q}\)^{- (q_{0,0}'-q_{0,0}'')/\ge}
\(W_R(x,\q',\gl)\)^\#
\(W_L(x,\q'',\gl)\)^\#\tau(x-\ge/2,\q'')
\right\}\frac{d\gl}{\gl}.
\een
Let $q_{0,0}'-q_{0,0}'' = m\ge$ and use that
\ben
&&
W_{L/R}(x,q_{0,0}'-m\ge, q_{0,1}'',\ldots,\gl) = 
W_{L/R}(x-m\ge,q_{0,0}', q_{0,1}'',\ldots,\gl).
\een
\qed



\end{document}